\documentclass [11pt]{article}
\usepackage{mathrsfs}
\usepackage{amsthm,amsmath,amsfonts,amssymb,graphicx,graphics,psfrag}

\newtheorem{thm}{Theorem}

\newtheorem{lem}{Lemma}
\newtheorem{ob}{Observation}
\newtheorem{cor}{Corollary}

\date{}

\begin{document}

\title{Uniformly Weighted Star-Factors of Graphs
\thanks{This work is supported by RFDP of Higher Education of China and Natural Sciences
and Engineering Research Council of Canada.}}

\author{ Yunjian Wu$^{1}$ and Qinglin Yu$^{1,2}$
\thanks{Email: yu@tru.ca (Q. Yu)}
\\ {\small $^1$ Center for Combinatorics, LPMC}
\\ {\small \small Nankai University, Tianjin, 300071, China}
\\ {\small  \small $^2$ Department of Mathematics and Statistics}
\\ {\small  Thompson Rivers University, Kamloops, BC, Canada}
}

\maketitle

\begin{abstract}
 A {\it star-factor} of a graph $G$ is a spanning subgraph of $G$
such that each component of which is a star. An {\it edge-weighting}
of $G$ is a function $w: E(G)\longrightarrow \mathbb{N}^+$, where
$\mathbb{N}^+$ is the set of positive integers. Let $\Omega$ be the
family of all graphs $G$ such that every star-factor of $G$ has the
same weights under a fixed edge-weighting $w$. In this paper, we
present a simple structural characterization of the graphs in
$\Omega$ that have girth at least five.
\begin{flushleft}
\end{flushleft}
\end{abstract}

\vskip 1mm

\vskip 1mm

\section{Introduction}
Throughout this paper, all graphs considered are simple. We refer
the reader to \cite{bb} for standard graph theoretic terms not
defined in this paper.

Let $G=(V,E)$ be a graph with vertex set $V(G)$ and edge set $E(G)$.
If $G$ is not a forest, then length of the shortest cycle in $G$ is
called the {\it girth} of $G$, denoted by $g(G)$ and a forest is
considered to have an infinite girth.   If $S \subset V(G)$, then
$G-S=G[V-S]$ is the subgraph of $G$ obtained by deleting the
vertices in $S$ and all the edges incident with them. Similarly, if
$E'\subset E(G)$, then $G-E'=(V(G), E(G)-E')$. We denote the degree
of a vertex $x$ in $G$ by $d_{G}(x)$, and the set of vertices
adjacent to $x$ in $G$ by $N_{G}(x)$. We also denote by $\delta(G)$
the minimum degree of vertices in $G$. A cycle $($or path$)$ with
$n$ vertices is denoted by $C_{n}$ $($or $P_{n})$. If vertices $u$
and $v$ are connected in $G$, the {\it distance} between $u$ and $v$
in $G$, denoted by $d_{G}(u,v)$, is the length of a shortest
$(u,v)$-path in $G$. The {\it diameter} of $G$ is the maximum
distance over all pairs of vertices in $G$. A {\it leaf} is a vertex
of degree one and a {\it stem} is a vertex which has at least one
leaf as its neighbor. A {\it star} is a tree isomorphic to $K_{1,n}$
for some $n \geq 1$, and the vertex of degree $n$ is called the {\it
center} of the star. A {\it star-factor} of a graph $G$ is a
spanning subgraph of $G$ such that each component of which is a
star. Clearly a graph with isolated vertices has no star-factors. On
the other hand, it is not hard to see that every graph without
isolated vertices admits a star-factor. If one limits the size of
the star used, the existence of such a star-factor is non-trivial.
In \cite{amahashi}, Amahashi and Kano presented a criterion for the
existence of a star-factor, i.e., \{$K_{1, 1}, \cdots , K_{1,
n}$\}-factor. Yu \cite{yu} obtained an upper bound on the maximum
number of edges in a graph with a unique star-factor.

An {\it edge-weighting} of a graph $G$ is a function $w:
E(G)\longrightarrow \mathbb{N}^+$, where $\mathbb{N}^+$ is the set
of positive integers. For a subgraph $H$, the {\it weight} of $H$
under $w$ is the sum of all the weight values for edges belonging to
$H$, i.e., $w(H)=\Sigma_{e\in E(H)}w(e)$. Motivated by the minimum
cost spanning trees and the optimal assignment problems, Hartnell
and Rall posed an interesting general question: for a given graph,
does there exist an edge-weighting function $w$ such that a certain
type of spanning subgraphs always have the same weights?  In
particular, they investigated the following narrow version of the
problem in
which the spanning subgraph is a star-factor.\\

\noindent\textbf{Star-Weighting Problem {{\rm{(Hartnell and Rall
\cite{hartnell}}}}}): For a given graph $G=(V, E)$, does there exist
an edge-weighting $w$ of $G$ such that every star-factor of $G$ has
the same weights under $w$?

\vspace{4mm}

  To start the investigation, one may consider a special case that
$w$ is a constant function, i.e., all edges in $G$ are assigned with
the same weights. In this case, two star-factors of $G$ have the
same weights if and only if they both have the same number of edges.
For simplicity, we assume that all edges are assigned with weight
one.

Let $\mathscr{U}$ be the family of all graphs $G$ such that if
$S_{1}$ and $S_{2}$ are any two star-factors of $G$, then $S_{1}$
and $S_{2}$ have the same number of edges. Clearly, $S_{1}$ and
$S_{2}$ have the same number of edges is equivalent to that they
have the same number of components. Hartnell and Rall
\cite{hartnell} classified the family $\mathscr{U}$ when graphs in
$\mathscr{U}$ have girth at least five. In \cite{wy}, the authors
characterized the family $\mathscr{U}$ when all its members have
girth three and minimum degree at least two.

We denote by $\Omega$ the family of all graphs $G$ such that every
star-factor of $G$ has the same weights under some fixed
edge-weighting $w$. In the definition of edge-weighting, we assume
that $w(e)\neq 0$ for every $e\in E(G)$. We also note that if $G\in
\mathscr{U}$ then $G\in \Omega$, but the converse is not always
true. In this paper, we give a structural characterization of the
graphs in $\Omega$ that have girth at least five.

\section{Main Result}
We start with a few easy observations and lemmas.

\vspace{3mm}

Note that if $H$ is a spanning subgraph of $G$, then any
star-factor of $H$ is also a star-factor of $G$.  The following
lemma will be used frequently in reducing the problem of
determining membership in $\Omega$ to its spanning subgraphs.

\begin{lem}\label{lem1}
Let $F$ be a subset of $E(G)$ such that $G-F$ has no isolated
vertices.  If $G-F \notin \Omega$, then $G \notin \Omega$ as well.
\end{lem}

    The above lemma implies that if $G$ is in $\Omega$, then so
is $G-F$. The basic idea to show that a graph does not belong to
$\Omega$ is to decompose $G$ into several components without
isolated vertices and then simply find one of them not belonging to
$\Omega$.

\begin{ob}\label{ob1}
Let $P_{6}=v_{1}v_{2}v_{3}v_{4}v_{5}v_{6}$ be a component of $G$
and $G\in \Omega$, then
$w(v_{3}v_{4})=w(v_{2}v_{3})+w(v_{4}v_{5})$.
\end{ob}

\begin{ob}\label{ob3}
Let $P_{5}=v_{1}v_{2}v_{3}v_{4}v_{5}$ be a component of $G$ and
$G\in \Omega$, then $w(v_{2}v_{3})=w(v_{3}v_{4})$.
\end{ob}

Use these two observations and Lemma \ref{lem1}, we get the
following observation.

\begin{ob}\label{ob2}
Let $C$ be a cycle.  If $|V(C)| = 6$ or $|V(C)| \geq 8$,  then $C
\notin \Omega$.
\end{ob}

We investigate the graphs, in  $\Omega$, with girth at least five
but without leaves first.

\begin{thm}\label{thm1}
If $\delta(G) \geq 2$ and either $g(G) = 6$ or $g(G) \geq 8$, then
$G\notin \Omega$.
\end{thm}
\noindent {\bf Proof.} Let $C$ be a cycle of order six in $G$ and
$F=\{u_{1}u_{2}\ | \ u_{1}\in V(C),u_{2}\in V(G)-V(C)\}$. Since the
girth of $G$ is six, there is no isolated vertices created in $G-F$.
So $C$, as a component in $G-F$, is a cycle of order six. Hence
$G\notin \Omega$ by Lemma \ref{lem1} and Observation
\ref{ob2}.

    For the case of $\delta(G) \geq 2$ and $g(G) \geq 8$, the argument
is very similar. \hfill$\Box$

\vspace{3mm}

\begin{thm}\label{thm3}
If $\delta(G) \geq 2$ and $g(G) = 7$, then $G\in \Omega$ if and only
if $G$ is a $7$-cycle.
\end{thm}
\noindent {\bf Proof.} Let $G$ be a $7$-cycle, let $w(e) = k$ for
each $e \in E(G)$ and $k\in \mathbb{N}^+$.  Then it is easy to check
that $G\in \Omega$.

    On the other hand, assume $G\in \Omega$ but $G$ is not a $7$-cycle. Let
$C = v_{1}v_{2}v_{3}v_{4}v_{5}v_{6}v_{7}$ be a cycle in $G$. Without
loss of generality, assume that $v_{7}$ has a neighbor $u$ not on
$C$. Let $F_{1}=\{u_{1}u_{2} \ | \ u_{1}\in V(C),u_{2}\in
V(G)-V(C)\}$. Since the girth of $G$ is seven, there are no isolated
vertices created in $G-F_{1}$. We see that $C$, as a component in
$G-F_{1}$, is a cycle of length seven. Since $G\in \Omega$, then all
edges on the cycle $C$ must have the same weights. Let
$F_{2}=\{u_{1}u_{2} \ | \ u_{1}\in
\{v_{1},v_{2},v_{3},v_{4},v_{5},v_{6}\},u_{2}\in
V(G)-\{v_{1},v_{2},v_{3},v_{4},v_{5},v_{6}\}\}$, then there are no
isolated vertices created in $G-F_{2}$ again. But
$P_{6}=v_{1}v_{2}v_{3}v_{4}v_{5}v_{6}$, as a component in $G-F_{2}$,
is a path of order six. Since $G\in \Omega$,
$w(v_{3}v_{4})=w(v_{2}v_{3})+w(v_{4}v_{5})$ by Observation
\ref{ob1}. Hence all the weights of edges on the cycle $C$ must be
$0$, a contradiction. \hfill$\Box$

\vspace{3mm}

\begin{lem}\label{lem2}
Let $G$ be a graph with an induced cycle of order five such that
four of the vertices are of degree two and the fifth is a stem. Then
$G \notin \Omega$.
\end{lem}

\noindent {\bf Proof.} Suppose $G$ belongs to $\Omega$, and $w$ is
an edge-weighting function such that all star-factors of $G$ have
the same weights under $w$. Let $C = v_{1}v_{2}v_{3}v_{4}v_{5}$ be a
$5$-cycle in $G$ with a stem $v_{5}$. Let $X$ be the set of leaves
adjacent to $v_{5}$ and $F_{1}=\{u_{1}v_{5} \ | \ u_{1}\in
V(G)-X-v_{4}\}$. A component $H_{1}$ of $G-F_{1}$ is isomorphic to
the graph shown in Figure $1(a)$. Since $G-F_{1}$ has no isolated
vertices, so $H_{1}\in \Omega$ and
$w(v_{3}v_{4})=w(v_{2}v_{3})+w(v_{4}v_{5})$ by Observation
\ref{ob1}. On the other hand, let $F_{2}=\{u_{1}v_{5} \ | \ u_{1}\in
V(G)-X-\{v_{1},v_{4}\}\}$, then a component $H_{2}$ of
$G-F_{2}-v_{1}v_{2}$ is isomorphic to the graph shown in Figure
$1(b)$. So $H_{2}\in \Omega$ and $w(v_{3}v_{4})=w(v_{4}v_{5})$ by
Observation \ref{ob3}. From the above two relations, we have
$w(v_{2}v_{3})=0$, a contradiction. \hfill$\Box$

\begin{figure}[h,t]
\setlength{\unitlength}{1pt}
\begin{center}
\begin{picture}(40,100)(0,60)

\put(-130,100){\line(1,0){120}} \put(-130,100){\line(-1,1){15}}
\put(-130,100){\line(-1,-1){15}}

\put(70,100){\line(1,0){90}}\put(70,100){\line(-1,1){15}}
\put(70,100){\line(-1,-1){15}}\put(70,100){\line(-1,0){18}}

\put(-130,90){$v_{5}$}\put(-100,90){$v_{4}$}
\put(-70,90){$v_{3}$}\put(-40,90){$v_{2}$} \put(-10,90){$v_{1}$}

\put(70,90){$v_{5}$}\put(100,90){$v_{4}$}
\put(130,90){$v_{3}$}\put(160,90){$v_{2}$} \put(45,90){$v_{1}$}

\put(-130,100){\circle* {4}}\put(-100,100){\circle* {4}}
\put(-70,100){\circle* {4}}\put(-40,100){\circle* {4}}
\put(-10,100){\circle* {4}} \put(-145,115){\circle* {4}}
\put(-145,85){\circle* {4}}

\put(70,100){\circle* {4}}\put(100,100){\circle* {4}}
\put(130,100){\circle* {4}}\put(160,100){\circle* {4}}
\put(53,100){\circle* {4}} \put(55,115){\circle* {4}}
\put(55,85){\circle* {4}}

\put(-90,70){$(a)$}\put(100,70){$(b)$}

\put(-5,50){Figure\ 1}

\end{picture}
\end{center}
\end{figure}

\vspace{3mm}

\begin{lem}\label{lem3}
Let $G$ be a graph in $\Omega$ with an induced $5$-cycle. If exactly
one of the vertices on this $5$-cycle has degree at least three,
then all of its neighbors not belonging to this $5$-cycle must be
stems.
\end{lem}

\noindent {\bf Proof.} Let $v$ be a vertex on the $5$-cycle of
degree at least three. Assume $v$ has a neighbor $x$ not on the
$5$-cycle and $x$ is not a stem. By Lemma \ref{lem2}, $x$ is not a
leaf. Let $F$ be the set of edges incident with $x$ except $vx$.
Then the graph $G-F$ has no isolated vertices, and the vertex $v$ is
a stem belonging to an induced $5$-cycle that satisfies the
hypothesis of Lemma \ref{lem2}. Thus $G \notin \Omega$, a
contradiction. \hfill$\Box$


\begin{thm}\label{thm4}
If $\delta(G) \geq 2$ and $g(G) = 5$,  then $G\in \Omega$ if and
only if $G$ is a $5$-cycle.
\end{thm}

\noindent {\bf Proof.}  If $G$ is a $5$-cycle, clearly $G \in
\Omega$ under a constant weight function.

    Next consider a graph $G \in \Omega$ with $\delta(G) \geq 2$
and $g(G) = 5$ but $G \ncong C_{5}$. Let $C =
v_{1}v_{2}v_{3}v_{4}v_{5}$ be a cycle in $G$. Assume, without loss
of generality, that $v_{5}$ has a neighbor $u$ not on $C$. Let
$F=\{u_{1}u_{2} \ | \ u_{1}\in
\{v_{1},v_{2},v_{3},v_{4},v_{5}\},u_{2}\in V(G)-V(C)-u\}$. If we
delete all edges in $F$ from $G$, then no isolated vertices created
in $G-F$ since $g(G) = 5$, so $G-F\in \Omega$ and $u$ is a stem in
$G-F$ by Lemma \ref{lem3}.

\vspace{0.5cm}

\begin{figure}[h,t]
\setlength{\unitlength}{1pt}
\begin{center}
\begin{picture}(40,100)

\put(-40,100){\line(-1,-1){20}} \put(-40,100){\line(1,-1){20}}
\put(-60,80){\line(0,-1){30}}\put(-20,80){\line(0,-1){30}}
\put(-60,51){\line(1,0){40}}

\put(-55,97){$v_{5}$} \put(-72,80){$v_{1}$} \put(-14,80){$v_{4}$}
\put(-72,50){$v_{2}$} \put(-14,50){$v_{3}$}

\put(-60,50){\circle* {4}} \put(-20,50){\circle* {4}}
\put(-60,80){\circle* {4}} \put(-40,100){\circle* {4}}
\put(-20,80){\circle* {4}}

\put(-60,50){\line(5,6){13}} \put(-40,80){\line(-1,-2){8}}

\put(-40,100){\line(0,-1){20}} \put(-40,80){\circle* {4}}
\put(-48,65){\circle* {4}}

\put(60,100){\line(-1,-1){20}} \put(60,100){\line(1,-1){20}}
\put(40,80){\line(0,-1){30}} \put(80,80){\line(0,-1){30}}
\put(40,51){\line(1,0){40}}

\put(45,97){$v_{5}$} \put(28,80){$v_{1}$} \put(86,80){$v_{4}$}
\put(28,50){$v_{2}$} \put(86,50){$v_{3}$}

\put(40,50){\circle* {4}} \put(80,50){\circle* {4}}
\put(40,80){\circle* {4}} \put(60,100){\circle* {4}}
\put(80,80){\circle* {4}}

\put(40,50){\line(5,6){13}}\put(68,65){\line(5,-6){13}}
\put(60,80){\line(-1,-2){8}}\put(60,80){\line(1,-2){8}}

\put(60,100){\line(0,-1){20}} \put(60,80){\circle* {4}}
\put(52,65){\circle* {4}}\put(68,65){\circle* {4}}

\put(-45,30){$(a)$}\put(55,30){$(b)$}

\put(-8,10){Figure\ 2}
\end{picture}
\end{center}
\end{figure}

\vspace{-0.5cm}

Moreover, $u$ has at most two leaves as its neighbors in $G-F$, and
they are adjacent to $v_{2}$ or (and) $v_{3}$ in $G$. Without loss
of generality, assume that all neighbors of $u$ except $v_{5}$ in
$G-F$ are leaves. Let $H$ be the component containing the cycle $C$
in $G-F$, and $H'$ be the induced subgraph of $G$ with $V(H)$. Then
$G-V(H')$ and $H'$ have no isolated vertices, and so $H'\in\Omega$
by Lemma \ref{lem1}. If there is exactly one leaf as a neighbor of
$u$ in $H$, then $H'$ is shown in Figure $2(a)$. Otherwise, $H'$ is
isomorphic to the graph shown in Figure $2(b)$. It is not hard to
check that both graphs are not in $\Omega$, a contradiction to $G\in
\Omega$. \hfill$\Box$

\vspace{3mm}

    From the four theorems above, we obtain the following corollary.
\begin{cor}\label{cor1}
If $G$ is a graph with $\delta(G) \geq 2$ and $g(G) \geq 5$, then
$G\in \Omega$ if and only if $G$ is a $5$-cycle or $7$-cycle.
Moreover, all edges of $G$ must have the same weights.
\end{cor}

    Next, we attempt to determine all members in $\Omega$ which have
girth at least five and with {\it leaves}. To derive our main
theorem, we need the following lemmas.

\begin{lem}\label{lem4}
Let $G$ be a graph of girth five and contain a $5$-cycle $C$ in
which no vertex is a stem and there exist two adjacent vertices of
degrees at least three. Then $G \notin \Omega$.
\end{lem}

\begin{figure}[h,t]
\setlength{\unitlength}{1pt}
\begin{center}
\begin{picture}(40,100)

\put(-40,100){\line(-1,-1){20}} \put(-40,100){\line(1,-1){20}}
\put(-60,80){\line(0,-1){30}}\put(-20,80){\line(0,-1){30}}
\put(-60,50){\line(1,0){40}}

\put(-60,51){\line(0,-1){20}} \put(-60,31){\line(-1,-1){10}}
\put(-60,31){\line(0,-1){10}} \put(-60,31){\line(1,-1){10}}
\put(-69,22){\circle* {4}} \put(-60,22){\circle* {4}}
\put(-51,22){\circle* {4}}

\put(-20,51){\line(0,-1){20}} \put(-20,31){\line(-1,-1){10}}
\put(-20,31){\line(1,-1){10}} \put(-29,22){\circle* {4}}
\put(-11,22){\circle* {4}} \put(-20,31){\circle* {4}}

\put(-75,31){$w_{1}$} \put(-15,31){$w_{2}$}

\put(-55,100){$v_{4}$} \put(-72,80){$v_{5}$} \put(-14,80){$v_{3}$}
\put(-72,50){$v_{1}$} \put(-14,50){$v_{2}$}

\put(-60,50){\circle* {4}} \put(-20,50){\circle* {4}}
\put(-60,80){\circle* {4}} \put(-40,100){\circle* {4}}
\put(-20,80){\circle* {4}} \put(-60,31){\circle* {4}}

\put(70,90){\line(-3,-1){30}}\put(70,90){\line(3,-1){30}}
\put(40,80){\line(0,-1){40}}\put(100,80){\line(0,-1){40}}
\put(40,41){\line(1,0){60}}\put(60,67){\line(1,0){20}}
\put(40,80){\line(3,-2){20}}\put(100,80){\line(-3,-2){20}}

\put(75,92){$v_{2}$}\put(28,80){$v_{3}$}\put(103,80){$v_{1}$}\put(28,35){$v_{4}$}\put(103,35){$v_{5}$}
\put(55,58){$y$}\put(80,58){$x$}

\put(70,90){\circle* {4}}\put(40,80){\circle* {4}}
\put(100,80){\circle* {4}}\put(40,40){\circle* {4}}
\put(100,40){\circle* {4}}\put(60,67){\circle* {4}}
\put(80,67){\circle* {4}}

\put(-45,5){$(a)$}\put(65,5){$(b)$}

\put(-8,-10){Figure\ 3}
\end{picture}
\end{center}
\end{figure}

\noindent {\bf Proof.} Assume $G \in \Omega$ and let $C=
v_{1}v_{2}v_{3}v_{4}v_{5}$ be the $5$-cycle in $G$ such that $v_{1}$
and $v_{2}$ both have degree at least three. Let $F_{1}$ be the set
of all edges not on $C$ but incident with one of $v_{3},v_{4}$ and
$v_{5}$. Since $g(G) = 5$, then no pair of vertices on $C$ have a
common neighbor not on $C$, and thus $G'=G-F_{1}$ has no isolated
vertices. Furthermore, neither $v_{1}$ nor $v_{2}$ is a stem in
$G'$. Let $G''$ be the graph obtained by deleting all edges incident
with $v_{1}$ but not on $C$ from $G'$, then none of stems created in
$G''$ is a neighbor of $v_{2}$. However, by Lemma \ref{lem1}, $G''$
belongs to $\Omega$ and so by Lemma \ref{lem3} all neighbors of
$v_{2}$ not on $C$, in $G''$, must be stems. Thus all neighbors of
$v_{2}$ not on $C$ in $G'$ must be stems. A similar argument yields
that all neighbors of $v_{1}$ not on $C$ in $G'$ must be stems. Let
$w_{1}$ and $w_{2}$ be stems adjacent to $v_{1}$ and $v_{2}$,
respectively, in $G'$. Then there exists at most one common neighbor
of degree two between $w_{1}$ and $w_{2}$ since $g(G)=5$. Let $X$ be
the set of leaves adjacent to $w_{1}$ or $w_{2}$ in $G'$ and
$F_{2}=\{u_{1}u_{2}\ | \ u_{1}\in
\{v_{1},v_{2},w_{1},w_{2}\},u_{2}\in
V(G')-\{v_{1},v_{2},w_{1},w_{2},v_{3},v_{5}\}-X\}$. Then $G'-F_{2}$
$($or $G'-F_{2}\cup uw_{1}$, if there exists a common neighbor $u$
of degree two between $w_{1}$ and $w_{2})$ has no isolated vertices
and a component $H$ in $G'-F_{2}$ $($or $G'-F_{2} \cup uw_{1})$ is
isomorphic to the graph shown in Figure $3(a)$. By Lemma \ref{lem1},
$H \in \Omega$. It is easy to show that $H\notin\Omega$ by
Observation \ref{ob1} and Observation \ref{ob3}, a
contradiction.\hfill$\Box$

\vspace{3mm}

\begin{lem}\label{lem5}
Let $G$ be a graph of girth five and contain a $5$-cycle $C=
v_{1}v_{2}v_{3}v_{4}v_{5}$ in which no vertex is a stem and two
nonadjacent vertices $v_{1}$ and $v_{3}$ of $C$ have degree at least
three. If $G \in \Omega$, then all neighbors of $v_{1}$ and $v_{3}$
not belonging to $C$ must be stems.
\end{lem}

\noindent {\bf Proof.} Assume that $v_{1}$ has a neighbor $x$ not on
$C$ but $x$ is not a stem. By Lemma \ref{lem4}, $v_{2},v_{4}$ and
$v_{5}$ must all have degree two in $G$. Let $F$ be the set of all
edges other than $xv_{1}$, that are incident with $x$, and let
$G'=G-F$. Since $x$ is not a stem in $G$, $G'\in \Omega$. If $v_{3}$
is not a stem in $G'$, let $F'$ be the set of edges, not in $C$, but
incident with $v_{3}$, then $G'-F' \in \Omega$ but it contains a
$5$-cycle satisfying the conditions of Lemma \ref{lem2}. Therefore
$v_{3}$ is a stem in $G'$. Moreover, we note that $v_{3}$ has
exactly one leaf, say $y$, as its neighbor in $G'$ since $g(G)=5$.
Thus, we arrive at an induced subgraph $H$ shown in Figure $3(b)$.
As $G-H$ has no isolated vertices, so $H\in \Omega$ by Lemma
\ref{lem1}. But by Theorem \ref{thm4}, $H \notin \Omega$, a
contradiction.\hfill$\Box$

\vspace{3mm}

\begin{lem}\label{lem6}
Let $G$ be a graph with an induced cycle of order six such that five
of the vertices are of degree two and the sixth is a stem. Then $G
\notin \Omega$.
\end{lem}

\noindent {\bf Proof.} Suppose $G \in \Omega$ under an
edge-weighting function $w$.  Let $C =
v_{1}v_{2}v_{3}v_{4}v_{5}v_{6}$ be a $6$-cycle in $G$ such that all
vertices on $C$ are of degree two except a stem $v_{1}$.

Let $X$ be the set of leaves adjacent to $v_{1}$ and
$F_{1}=\{u_{1}v_{1} \ | \ u_{1}\in V(G)-X-\{v_{2},v_{6}\}\}$. A
component $H_{1}$ of $G-F_{1}-v_{5}v_{6}$ is isomorphic to the graph
shown in Figure $4(a)$. Since $G-F_{1}-v_{5}v_{6}$ has no isolated
vertices, so $H_{1}\in \Omega$ and
$w(v_{2}v_{3})=w(v_{1}v_{2})+w(v_{3}v_{4})$ by Observation
\ref{ob1}. On the other hand, let $F_{2}=\{u_{1}v_{1} \ | \ u_{1}\in
V(G)-X-v_{2}\}$. A component $H_{2}$ of $G-F_{2}-v_{4}v_{5}$ is
isomorphic to the graph shown in Figure $4(b)$. So $H_{2}\in \Omega$
and $w(v_{1}v_{2})=w(v_{2}v_{3})$ by Observation \ref{ob3}. From the
above two relations, we have $w(v_{3}v_{4})=0$, a contradiction.
\hfill$\Box$

\vspace{-1cm}
\begin{figure}[h,t]
\setlength{\unitlength}{1pt}
\begin{center}
\begin{picture}(40,100)(0,60)

\put(-130,100){\line(1,0){120}} \put(-130,100){\line(-1,1){15}}
\put(-130,100){\line(-1,-1){15}}\put(-130,100){\line(-1,0){15}}

\put(70,100){\line(1,0){90}}\put(70,100){\line(-1,1){15}}
\put(70,100){\line(-1,-1){15}}

\put(-150,90){$v_{6}$} \put(-130,90){$v_{1}$}\put(-100,90){$v_{2}$}
\put(-70,90){$v_{3}$}\put(-40,90){$v_{4}$} \put(-10,90){$v_{5}$}

\put(70,90){$v_{1}$}\put(100,90){$v_{2}$}
\put(130,90){$v_{3}$}\put(160,90){$v_{4}$}

\put(-130,100){\circle* {4}}\put(-100,100){\circle* {4}}
\put(-70,100){\circle* {4}}\put(-40,100){\circle* {4}}
\put(-10,100){\circle* {4}} \put(-145,115){\circle* {4}}
\put(-145,85){\circle* {4}}\put(-145,100){\circle* {4}}

\put(70,100){\circle* {4}}\put(100,100){\circle* {4}}
\put(130,100){\circle* {4}}\put(160,100){\circle* {4}}
\put(55,115){\circle* {4}} \put(55,85){\circle* {4}}

\put(-90,70){$(a)$}\put(100,70){$(b)$}

\put(-5,50){Figure\ 4}

\end{picture}
\end{center}
\end{figure}

We now characterize the graphs in $\Omega$ of gith at least five.
The results above reduce the problem to considering those graphs
that have at least one leaf. The following lemma is true regardless
of girth.

\begin{lem}\label{lem7}
If $G$ is a graph in which every vertex is either a leaf or a stem,
then $G$ belongs to $\Omega$
\end{lem}

\noindent {\bf Proof.} Let $L$ be the set of leaves of $G$ and $W$
be the set of stems. Any star-factor of $G$ must contain all edges
joining a vertex in $L$ and a vertex in $W$ but cannot contain any
edge incident with two vertices of $W$. Hence $G$ has a unique
star-factor and $G\in \Omega$.\hfill$\Box$

\vspace{3mm}

    Now we prove our main result.

\begin{thm}\label{thm5}
Let $G$ be a connected graph of girth at least five. Then $G\in
\Omega$ if and only if $G$ is

    $1)$ a $C_5$, or

    $2)$ a $C_7$, or

    $3)$ $G$ has leaves and each vertex in $G$ is either a leaf or a stem, or

    $4)$ $G$ has leaves but contains at least one vertex which is neither a leaf nor a stem,
then each component of the graph obtained by removing the leaves and
stems from $G$ is one of the following:

 $4a)$ a $5$-cycle with at most two vertices of degree three or more
in $G$. Furthermore, if there are two such vertices, then they are
non-adjacent on the $5$-cycle;

 $4b)$ a star $K_{1,m}$ $(m\geq 1)$. Moreover, the center of
$K_{1,m}$ has degree $m$ in $G$ for $m\geq 2$;

 $4c)$ an isolated vertex.
\end{thm}

\noindent {\bf Proof.} If $g(G) \geq 5$ and $\delta(G)\geq 2$, then
the theorem follows from Corollary \ref{cor1}.

    If $\delta(G) = 1$ and each vertex in $G$ is either a leaf or a stem, then
$G\in \Omega$ by Lemma \ref{lem7}.

     Next, suppose $\delta(G) = 1$ {\it and} $G\in
\Omega$, but $G$ contains at least one vertex which is neither a
leaf nor a stem. Let $L$ be the set of leaves and $S=N(L)$ the set
of stems. If $G-(L\cup S)$ has a component which is a 5-cycle $C$,
then, by Lemmas \ref{lem4} and \ref{lem5}, there are at most two
vertices of $C$ with degree greater than two in $G$ $($and if there
are two such vertices, they are nonadjacent on $C)$ and all of their
neighbors not on $C$ must be stems.

Hence we now consider those components of $G-(L\cup S)$ that have no
$5$-cycles. Let $H$ be such a component, then $g(H)\geq 6$. We shall
show that {\it the diameter of $H$ is at most two}.

Suppose the diameter of $H$ is at least three, then there exists a
path $P = abcd$ in $H$ such that $a$ is adjacent to a stem $s$ of
$G$.

\begin{figure}[h,t]
\setlength{\unitlength}{1pt}
\begin{center}
\begin{picture}(40,100)

\put(50,80){\line(1,0){90}}\put(50,80){\line(0,-1){20}}
\put(50,60){\line(-1,-2){10}}\put(50,60){\line(1,-2){10}}

\put(48,85){$a$}\put(78,85){$b$}\put(108,85){$c$}\put(138,85){$d$}\put(43,60){$s$}

\put(50,80){\circle* {4}}\put(80,80){\circle* {4}}
\put(110,80){\circle* {4}}\put(140,80){\circle* {4}}
\put(50,60){\circle* {4}}\put(60,40){\circle* {4}}
\put(40,40){\circle* {4}}

\put(-90,80){\line(1,0){60}}\put(-90,80){\line(0,-1){20}}
\put(-90,60){\line(-1,-2){10}}\put(-90,60){\line(1,-2){10}}

\put(-92,85){$a$}\put(-62,85){$b$}\put(-32,85){$c$}\put(-97,60){$s$}

\put(-90,80){\circle* {4}}\put(-60,80){\circle* {4}}
\put(-30,80){\circle* {4}}\put(-100,40){\circle* {4}}
\put(-90,60){\circle*{4}}\put(-80,40){\circle* {4}}

\put(-70,10){$(a)$}\put(80,10){$(b)$}

\put(-5,-10){Figure\ 5}
\end{picture}
\end{center}
\end{figure}

\vspace{2mm}

{\it Claim 1.} There exist no common neighbors of degree two between
$a$ and $d$, or $s$ and $c$, or $s$ and $d$ in $G$.

Since $g(H)\geq 6$, there is no common neighbor of degree two
between $a$ and $d$.

Let $X$ be the set of leaves adjacent to vertices $s$ in $G$.
Suppose there is a common neighbor $u$ of degree two between $s$ and
$c$. Let $F_{1}=\{u_{1}u_{2} \ | \ u_{1}\in \{a,b,c,u,s\},u_{2}\in
\{V(G)-\{a,b,c,u,s\}-X\}\}$. Then the graph $G_{1}=G-F_{1}$ has no
isolated vertices and has a component satisfying the hypothesis of
Lemma \ref{lem2}, so $G \notin \Omega$, a contradiction.

Suppose $v$ is a common neighbor of degree two between $s$ and $d$.
Let $F_{2}=\{u_{1}u_{2} \ | \ u_{1}\in \{a,b,c,d,s,v\},u_{2}\in
\{V(G)-\{a,b,c,d,s,v\}-X\}\}$. Then $G_{2}=G-F_{2}$ has no isolated
vertices but has a component satisfying the hypothesis of Lemma
\ref{lem6}, a contradiction to $G\in \Omega$ by Lemma \ref{lem1}.

\vspace{2mm}

    Claim $1$ yields that there are no common neighbors
of degree two between any two vertices of $\{s,a,b,c,d\}$. Let
$F_{3}=\{u_{1}u_{2} \ | \ u_{1}\in \{a,b,c,s\},u_{2}\in
\{V(G)-\{a,b,c,s\}-X\}\}$, Since $g(G)\geq 5$, then $G_{3}=G-F_{3}$
has no isolated vertices and has a component $H'$ isomorphic to the
graph shown in Figure $5(a)$. By Observation \ref{ob3}, we have

\vspace{-3mm}

\begin{equation}
w(ab)=w(as),
\end{equation}

On the other hand, let $F_{4}=\{u_{1}u_{2} \ | \ u_{1}\in
\{a,b,c,d,s\},u_{2}\in \{V(G)-\{a,b,c,d,s\}-X\}\}$, then
$G_{4}=G-F_{4}$ has no isolated vertices and has a component $H''$
isomorphic to the graph shown in Figure $5(b)$. Then $H''\in \Omega$
and

\vspace{-3mm}

\begin{equation}
w(ab)=w(as)+w(bc).
\end{equation}

Equations $(1)$ and $(2)$ imply that $w(bc)=0$, a contradiction.\\

Hence the diameter of $H$ is at most two, i.e., $H$ is either an
isolated vertex or isomorphic to a star, say $K_{1,m}$. For $m\geq
2$, let the vertices of $H$ be $c,b_{1},b_{2},\cdots, b_{m}$ where
$c$ has degree $m$ in $H$. For each $1\leq i\leq m$, let $s_{i}$ be
a stem of $G$ adjacent to $b_{i}$.

\vspace{2mm}

    {\it Claim 2.} If $m \geq 2$, then $c$ does not have a neighbor, in $G$, which is a stem.

    Otherwise, let $s$ be one of such neighbors and let $L_{i}$ and $L_s$ be
the sets of leaves adjacent to $s_{i}$ and $s$ in $G$, respectively.
If there exists a vertex $u$ adjacent only to vertices in $\{s_{1},
s_2, \dots, s_m, b_1, b_2, \dots, b_m\}$ and to at least one vertex
of $\{s_{1}, s_2, \dots, s_m\}$, then we can delete all edges which
are adjacent to $u$ except one $s_{k}u$ $($for some $1\leq k\leq
m)$. Thus we obtain a spanning subgraph of $G$ without isolated
vertices, and $u$ is a leaf adjacent to $s_{k}$. Hence we may assume
that there are no vertices only adjacent to vertices $\{s_{1}, s_2,
\dots, s_m, b_1, b_2, \dots, b_m\}$. For the same reason, we may
assume no vertices only adjacent to vertices $s$ and $s_{i},1\leq
i\leq m$. Let $F_{5}=\{u_{1}u_{2} \ | \ u_{1}\in
\{c,b_{1},b_{2},\cdots, b_{m},s_{1},s_{2},\cdots,s_{m},s\},u_{2}\in
\{V(G)-\{c,b_{1},b_{2},\cdots,
b_{m},s_{1},s_{2},\cdots,s_{m},s\}-L_{1}-L_{2}-\cdots-L_{m}-L_s\}\}$.
Then $G_{5}=G-F_{5}$ has no isolated vertices, since $g(G)\geq 5$,
but has a component $H'''$ isomorphic to the graph shown in Figure
$6$. Let the total weights of all edges incident with the leaves in
$L_{1}\cup L_{2} \cup \cdots \cup L_{m} \cup L_s$ be $w'$. Since
$H'''\in \Omega$, then by Observation \ref{ob1}, we have

\vspace{-3mm}

\begin{equation}
w(cb_{i})=w(b_{i}s_{i})+w(cs), \ 1\leq i\leq m.
\end{equation}

Since $G\in \Omega$, we also have

\vspace{-6mm}

\begin{equation}
w'+w(cb_{1})+w(cb_{2})+\cdots+w(cb_{m})=w'+w(cs)+w(b_{1}s_{1})+w(b_{2}s_{2})+\cdots+w(b_{m}s_{m}).
\end{equation}

    Equations $(3)$ and $(4)$ imply $m=1$, a contradiction to
$m\geq 2$. \\

    From Claim 2, we conclude that the center $c$ of the star $H=K_{1,m}\ (m\geq 2)$ has degree $m$
in $G$.

\begin{figure}[h,t]
\setlength{\unitlength}{1pt}
\begin{center}
\begin{picture}(40,100)

\put(20,50){\line(-1,1){15}}

\put(20,50){\line(-1,-1){15}}\put(5,65){\line(-1,0){15}}
\put(5,65){\line(-1,1){10}} \put(5,35){\line(-1,-1){10}}
\put(-5,26){\line(-2,-1){15}}\put(-5,26){\line(-1,-2){8}}

\put(20,50){\line(0,-1){15}}\put(20,35){\line(0,-1){15}}
\put(20,20){\line(-1,-1){10}}\put(20,20){\line(1,-1){10}}

\put(20,50){\line(1,1){15}}\put(35,65){\line(1,1){15}}
\put(50,80){\line(0,1){15}}\put(50,80){\line(1,1){15}}
\put(50,80){\line(1,0){15}}

\put(23,48){$c$}\put(8,68){$s$}\put(-3,40){$b_{1}$}\put(-12,32){$s_{1}$}\put(24,33){$b_{2}$}\put(22,22){$s_{2}$}
\put(40,60){$b_{m}$}\put(50,72){$s_{m}$}

\put(20,50){\circle* {4}}\put(5,65){\circle*{4}}
\put(-10,65){\circle* {4}}\put(-5,75){\circle* {4}}
\put(-5,26){\circle* {4}}\put(-20,19){\circle* {4}}
\put(-13,10){\circle* {4}}\put(20,35){\circle* {4}}
\put(20,20){\circle* {4}}\put(10,10){\circle* {4}}
\put(30,10){\circle* {4}}\put(35,65){\circle* {4}}
\put(50,80){\circle* {4}}\put(65,80){\circle*{4}}
\put(65,95){\circle*{4}}\put(50,95){\circle* {4}}
\put(5,35){\circle* {4}}\put(45,30){\circle* {2}}
\put(52,35){\circle* {2}}\put(55,45){\circle* {2}}

\put(10,-20){Figure\ 6}
\end{picture}
\end{center}
\end{figure}

\vspace{0.5cm}

Therefore every
component of $G-(L\cup S)$ is one of 4a), 4b) or 4c). \\

Conversely, assume $G$ has the specified structure. In the
following, we present an edge-weighting function such that every
star-factor of $G$ has the same weights.

\vspace{2mm}

{\it Case 1.} No component of $G-(L\cup S)$ is $K_{1,1}$.

In this case, all edges of $G$ are assigned the same weight. We only
need to show that all star-factors of $G$ have the same number of
edges. Let $T$ be any star-factor of $G$. Then $T$ contains exactly
one edge incident to one leaf of $G$. Let $H$ be a component of
$G-(L\cup S)$. If $H$ is a $5$-cycle, then $T$ either contains three
edges of $H$ or two edges of $H$ and exactly one edge of $T$ joining
$H$ to a stem of $G$. If $H=K_{1,m}$ ($m\geq 2$), then $T$ contains
precisely $m$ edges incident with a vertex of $H$. In particular,
for each leaf $x$ of $H$ either the edge joining $x$ to the center
of $H$ or exactly one edge joining $x$ to a stem of $G$ must be in
$T$. Also note that at least one edge of $H$ must be in $T$.
Finally, if $H$ is an isolated vertex $u$, then $T$ contains exactly
one edge joining $u$ to a stem of $G$.

\vspace{2mm}

{\it Case 2.} $K_{1,1}$ appears as a components of $G-(L\cup S)$.

    For each $K_{1,1}=uv$ of $G-(L\cup S)$, assign edge-weights
for the edges adjacent to $N_G(u) \cup N_G(v)$ as follows:

$$w(e)=\left\{\begin{array}{lllcrrr}
    a & e\in \{ux \ | \ x \in N(u)\}, e \not = uv\\
    b  & e\in \{vx \ | \ x \in N(v)\}, e \not = uv\\
    a+b& e=uv
\end{array} \right.$$

\noindent where $a, b >0$. All other edges are assigned the same
weight.

     Let $T$ be any star-factor of $G$. Then $T$
contains exactly one edge incident to one leaf of $G$. For each
component which is a $5$-cycle or an isolated vertex, it can be
dealt with as in Case $1$. For $H=uv$, then $T$ contains an edge
$uv$ or an edge joining a stem $s_{1}$ to $u$ and an edge joining
another stem $s_{2}$ to $v$. Since $w(uv)=w(us_{1})+w(vs_{2})$, we
conclude that every star-factor of $G$ has the same weights.

This completes the proof.\hfill$\Box$

\vspace{5mm}

 From the proof of theorem above, we obtain the following corollary.
\begin{cor}\label{cor2}
If $G$ is a tree, then $G\in \Omega$ if and only if each component
of the graph obtained by removing the leaves and stems from $G$ is
empty or a star $K_{1,m}$ $(m\geq 1)$ with center having degree $m$
in $G$ for $m\geq 2$ or an isolated vertex.
\end{cor}

\vspace{2mm}

\noindent {\bf Remark 1.} The graph $G$ shown in Figure 7 is an
example which is in $\Omega$ but requiring non-constant edge-weight
function. To see this, assume that all the edges have the same
weights, then we can find two star-factors with $10$ edges and $7$
edges, respectively. So $G \notin \mathscr{U}$. But if we give a
non-constant edge-weight function $w$ as follows:

$$w(e)=\left\{\begin{array}{llcrr}
    2k & e\in \{a,b,c\}\\
    k  & e\in \{E(G)-\{a,b,c\}\}\\
\end{array} \right.$$
\noindent where $k>0$. It is not hard to verify that all the
star-factors of $G$ have the same weights under $w$.

\begin{figure}[h,t]
\setlength{\unitlength}{1pt}
\begin{center}
\begin{picture}(40,100)

\put(15,50){\line(-1,1){15}}\put(15,50){\line(-1,-1){15}}
\put(15,50){\line(1,0){25}}
\put(40,50){\line(1,1){15}}\put(40,50){\line(1,-1){15}}
\put(0,65){\line(0,1){15}}\put(0,35){\line(0,-1){15}}
\put(0,66){\line(-1,0){25}}\put(0,35){\line(-1,0){25}}
\put(-25,65){\line(0,-1){30}}

\put(55,65){\line(0,1){15}}\put(55,35){\line(0,-1){15}}
\put(55,66){\line(1,0){25}}\put(55,35){\line(1,0){25}}
\put(80,65){\line(0,-1){30}}

\put(15,50){\circle* {4}}\put(40,50){\circle* {4}}
\put(55,65){\circle* {4}}\put(55,35){\circle* {4}}
\put(0,65){\circle* {4}}\put(0,35){\circle* {4}}
\put(0,80){\circle*{4}}\put(0,20){\circle* {4}}
\put(-25,65){\circle*{4}}\put(-25,35){\circle* {4}}

\put(55,80){\circle* {4}}\put(55,20){\circle* {4}}
\put(80,65){\circle*{4}}\put(80,35){\circle* {4}}

\put(-32,46){$a$}\put(25,40){$b$}\put(83,46){$c$}

\put(10,-10){Figure\ 7}
\end{picture}
\end{center}
\end{figure}

\vspace{2mm}

\noindent {\bf Remark 2.} The main theorem has classified all graphs
in $\Omega$ with girth at least five. The families remaining to be
determined are graphs of girth three or four. It seems that the
structures of both families are much more complicated, but it would
be an interesting problem to investigate.

%

\end{document}